\documentclass[10pt]{amsart}
\pdfoutput=1
\usepackage{amsmath, amsthm, amssymb,stmaryrd}
\usepackage{ifpdf}
\usepackage[pdftex]{graphicx}
\usepackage{tikz}
\usetikzlibrary{matrix,arrows,calc}
\usepackage[all]{xy}

\usepackage[pdftex,plainpages=false,hypertexnames=false,pdfpagelabels]{hyperref}
 \theoremstyle{plain}
 \newtheorem{thm}{Theorem}[section]

 \newtheorem{lem}[thm]{Lemma}
 \newtheorem{prop}[thm]{Proposition}
 \theoremstyle{definition}
 \newtheorem{defn}[thm]{Definition}
 
 \newtheorem{ex}[thm]{Example}
 \theoremstyle{remark}
 \newtheorem{rmk}[thm]{Remark}

\def\beq{\begin{eqnarray}}
\def\eeq{\end{eqnarray}}

\newcommand{\be}{\begin{equation}}
\newcommand{\ee}{\end{equation}}

 \newcommand{\bp}{\begin{proof}[Proof]}
 \newcommand{\ep}{\end{proof}}
\DeclareMathOperator{\SM}{\underline{\sf SMfld}}

\DeclareMathOperator{\ev}{{\rm ev}}
\DeclareMathOperator{\shol}{{\rm sHol}}

\newcommand{\sq}{\mathord{/\!\!/}}

\def\tt{{\sf {t}}}
\def\R{{\mathbb{R}}}

\def\E{{\mathbb{E}}}
\def\C{{\mathbb{C}}}

\def\Z{{\mathbb{Z}}}

\def\sPath{ {{\sf sP}}}
\def\sLoop{ {{\sf sL}}}

\def\path{ {{\mathcal P}_0}}

\def\pt{{\rm pt}}

\def\Ad{{\rm Ad}}
\def\id{{\rm id}}
\def\End{\mathop{\sf End}}
\def\Hom{\mathop{\sf Hom}}

\begin{document}

\title{The equivariant Chern character as super holonomy on loop stacks}

\author{Daniel Berwick-Evans and Fei Han}

\address{Department of Mathematics, University of Illinois at Urbana--Champaign and Department of Mathematics, National University Singapore}

\email{danbe@illinois.edu and mathanf@nus.edu.sg}

\date{\today}

\begin{abstract}
We study super parallel transport around super loops in a quotient stack, and show that this geometry constructs a global version of the equivariant Chern character. 
\end{abstract}

\maketitle

\section{Introduction}

Let $V\to M$ be an equivariant vector bundle. The image of the equivariant K-theory class $[V]\in {\rm K}_G(M)$ under complexification, 
\beq
{\rm K}_G(M)\to {\rm K}_G(M)\otimes \C\label{eq:equivChern}
\eeq
is a version of the equivariant Chern character, called the \emph{delocalized} or \emph{global equivariant} Chern character, e.g., see~\cite{Rosu,FHT}. By the Atiyah--Segal completion theorem~\cite{AS}, it encodes strictly more information than the equivariant Chern character valued in 2-periodic Borel equivariant cohomology, ${\rm H}_G(M;\C)$. Indeed, when $M=\pt$ the map~\eqref{eq:equivChern} takes the character of (virtual) $G$-representations as functions on~$G$, and the Taylor expansion of this function at the identity~$e\in G$ is the Chern character valued in~${\rm H}_G(\pt;\C)$. More generally, for each $g\in G$, the image of~$[V]$ under~\eqref{eq:equivChern} remembers information about the restriction of~$V$ to each fixed point set~$M^g$ analogous to the information of the character of a representation evaluated at~$g$. 

Endowing $V$ with an invariant connection~$\nabla$, Duflo--Vergne~\cite{DV} and Block--Getzler \cite{BlockGetzler} construct a differential-geometric description of~\eqref{eq:equivChern}. Owing to the Chern characters obtained at each $g\in G$, they call this a \emph{bouquet of Chern characters}. Among other applications, they explain how this plays an important role in certain versions of geometric quantization~\cite{DV, Vergne}. Their construction (see~\S\ref{appen}) uses Chern--Weil theory, starting with the definition of the \emph{equivariant curvature form}. Typically this definition is motivated by desirable algebraic features in connection with the algebra of equivariant differential forms. Although this approach has beautiful consequences, it makes the geometric origin of the bouquet of Chern characters somewhat obscure. 

One geometric interpretation of (ordinary) curvature comes from parallel transport around infinitesimal loops. The goal of this paper is to show that for a manifold~$M$ with a $G$-action, the bouquet of Chern characters (and with it, the equivariant curvature) comes from parallel transport around an appropriate version of infinitesimal super loops in the quotient stack~$M\sq G$. From this perspective, the complexity of these objects is a reflection of the fact that the loop space of $M\sq G$ is richer than the loop space of $M$ with a $G$-action. An \emph{equivariant loop} has as data a path $\gamma\colon [0,1]\to M$ together with $g\in G$ such that $\gamma(0)=g\cdot \gamma(1)$. In particular, \emph{constant} equivariant loops encode all the fixed point sets~$M^g$ for all $g\in G$. We use super geometry to define a precise notion of infinitesimal loops, so that the bouquet of Chern characters is a function on a certain equivariant \emph{super} loop space. The analogous result in the non-equivariant setting is the construction of the (ordinary) Chern character via super holonomy in the second named author's thesis~\cite{Han}. 

Let~$\mu\colon G\times M\to M$ be the action of a compact Lie group~$G$ on a smooth manifold~$M$ and~$(V,\mu^V, \nabla)$ be an equivariant vector bundle with action map $\mu^V\colon G\times V\to V$ and invariant connection $\nabla$.  Let $p: G\times M \to M$ be the projection, and $\mathcal{Z}(g)$ denote the Lie algebra of the stabilizer of $g$.

\begin{thm} 
Equivariant super  parallel transport defines a section
$$
{\rm sPar}_G(V,\nabla)\in \Omega^\bullet(G\times\mathfrak{g}\times M;{\sf Hom}(p^*V,\mu^*V))
$$
that when restricted to $\{g\}\times \mathcal{Z}(g)\times M^g$ determines a section of the endomorpism bundle of $p^*V$ whose trace is the contribution to the bouquet of  Chern characters of~$(V,\mu^V, \nabla)$ at~$g\in G$.
These traces are the equivariant super holonomy of~$V$ on constant equivariant super loops in~$M$.\label{thm1}
\end{thm}

Our conceptual framework is motivated by Stolz and Teichner's program to realize supersymmetric field theories as cocycles for generalized cohomology theories. Just as geometric models for de~Rham cohomology and K-theory lead to interesting equivariant refinements, one intriguing aspect of Stolz and Teichner's proposal is a connection between \emph{gauged} field theories and \emph{equivariant} cohomology theories (\cite{ST11, HSST}). Theorem~\ref{thm1} can be understood as a concrete verification of this idea in one example: the bouquet of Chern characters determines a class in ${\rm K}_G(M)\otimes \C$, and equivariant super paths are the fields for the version of gauged super symmetric quantum mechanics relevant to the physics proof of the equivariant index theorem~\cite{Alvarez}. The ordinary Chern character has a loop space lifting called the Bismut--Chern character \cite{B85}, which can be interpreted as a dimension reduction functor for supersymmetric Euclidean field theories \cite{HST}. We will study the equivariant  Bismut--Chern character \cite{B85} in the framework of gauged field theories in a forthcoming project. 

$\, $

\noindent {\bf Acknowlegement} We are indebted to Prof. Stephan Stolz and Prof. Peter Teichner for their insight and support. F. H. would also like to thank Prof. Varghese Mathai and Prof. Weiping Zhang for helpful discussions, and the support from National University of Singapore during his sabbatical leave as well as the hospitality of Max-Planck Institute for Mathematics at Bonn and Chern Institute at Tianjin during his visits. 

\section{Super parallel transport and the (ordinary) Chern character}

We start by reviewing the construction of the Chern character via super parallel transport in a manner that generalizes to the equivariant setting. We follow the conventions for supermanifolds and super Lie groups in~\cite{DM}. Super parallel transport was first constructed by Dumitrescu~\cite{Florin}. 

\subsection{The super Lie group $\E^{1|1}$}

Let $\E^{1|1}$ denote the super Lie group whose underlying super manifold is $\R^{1|1}$ with group structure
$$
(t,\theta)\cdot (s,\eta)=(t+s+\theta\eta,\theta+\eta),\quad (t,\theta),(s,\eta)\in \R^{1|1}(S). 
$$
The Lie algebra of $\E^{1|1}$ is freely generated (as a super Lie algebra) by a single odd element $D=\partial_\theta+\theta\partial_t$. The reduced group of $\E^{1|1}$ is $\E\cong \R$ with the usual group operation from addition. The Lie algebra of $\E$ is free on a single \emph{even} generator, so in this sense $\E^{1|1}$ provides an odd or super-analog of this familiar Lie group. 

By differentiating at the identity, an $\E^{1|1}$-action on a super manifold $S$ determines an odd vector field that is the infinitesimal action by~$D$. Conversely, given an odd vector field~$X$ on a super manifold~$S$ for which the (even) vector field $X^2=\frac{1}{2}[X,X]$ on the reduced manifold of~$S$ generates a flow, we obtain an $\E^{1|1}$-action on $S$. This follows from the description of super Lie group actions in terms of the action of the reduced group (in this case~$\R$) and a compatible action by the super Lie algebra of the group~\cite{DM}. 

\begin{ex}
A prototypical $\E^{1|1}$-action is the left action~$\tt\colon \E^{1|1}\times \R^{1|1}\to \R^{1|1}$, which we call the \emph{super translation action}. The odd vector field associated with this action is~$D$. We call the action by $\E<\E^{1|1}$ on $\R^{1|1}$ the \emph{even translation} action. 
\end{ex}

\begin{ex} For a compact manifold~$M$ with $G$-action, any element $X\in \mathfrak{g}$ defines an $\E^{1|1}$-action on $\Pi TM$ generated by the odd derivation $d-\iota_X$ acting on $\Omega^\bullet(M)\cong C^\infty(\Pi TM)$. 
\end{ex}

\subsection{Super path and super loop spaces}\label{sec:susypath}

Define the \emph{super path space} $\sPath(M)$ as the presheaf on super manifolds
$$
\sPath(M):=\SM(\R^{1|1},M)=\{ S\mapsto {\sf SMfld}(\R^{1|1}\times S,M)\}.
$$
The \emph{constant super paths} are the subsheaf of $\sPath(M)$ invariant under even translations so that they factor as
$$
\gamma\colon \R^{1|1}\times S\twoheadrightarrow (\R^{1|1}\times S)/\E\cong \R^{0|1}\times S\to M,
$$
where the first arrow is the fiberwise quotient map by the action of $\E<\E^{1|1}$ on $\R^{1|1}$. The universal family for these constant maps is the super manifold $S=\SM(\R^{0|1},M)$, and the universal map $\gamma$ is determined by the evaluation map
\beq
\ev\colon \R^{0|1}\times \SM(\R^{0|1},M)\to M,\label{eq:ev}
\eeq
adjoint to the identity map. On functions, after identifying $C^\infty(\SM(\R^{0|1},M))\cong C^\infty(\Pi TM)\cong \Omega^\bullet(M)$ and $C^\infty(\R^{0|1})\cong \R[\theta]$, we have $\ev^*(f)=f+\theta df$. 

An $S$-point $\gamma$ of $\sPath(M)$ defines an \emph{$S$-family of super loops} if the diagram commutes
\beq
\begin{tikzpicture}[baseline=(basepoint)];
\node (A) at (0,0) {$\Z\times \R^{1|1}\times S$};
\node (B) at (4,-0.75) {$\Z\times M$};
\node (C) at (0,-1.5) {$ \R^{1|1}\times S$};
\draw[->] (A) to node [above=1pt]{$\gamma\circ {\rm pr}$} (B);
\draw[->] (A) to node [left=1pt] {${\sf t}_\Z$} (C);
\draw[->] (C) to node [below=1pt] {$\gamma$} (B);
\path (0,-.75) coordinate (basepoint);
\end{tikzpicture}\nonumber
\eeq
where ${\sf t}_\Z$ denotes the action by translations of $\Z\subset \E\subset \E^{1|1}$ on $\R^{1|1}$. We denote this sub-presheaf by $\sLoop(M)\subset \sPath(M)$ and call it the \emph{super loop space} of~$M$. The inclusion of the constant super paths $\SM(\R^{0|1},M)$ into $\sPath(M)$ factors through the super loops, 
$$
\SM(\R^{0|1},M)\hookrightarrow \sLoop(M)\subset \sPath(M). 
$$

\subsection{Super parallel transport}

Let $(V,\nabla)$ be a vector bundle on $M$, and $\gamma\colon \R^{1|1}\times S\to M$ a family of super paths. 

\begin{defn}
A section $s\in \Gamma(\R^{1|1}\times S,\gamma^*V)$ is \emph{super parallel} if 
$$
(\gamma^*\nabla)_Ds=0
$$
where $D=\partial_\theta+\theta\partial_t$.
\end{defn}

To get a better flavor for the super parallel transport equation, it is instructive to study it in \emph{components}, which are defined as
\beq
s_0:=i_0^*s\in \Gamma(\R\times S,i_0^*\gamma^*V),\quad s_1:=i_0^* ((\gamma^*\nabla)_Ds)\in \Gamma(\R\times S,i_0^*\gamma^*\Pi V)\label{eq:components}
\eeq
for $s\in \Gamma(\R^{1|1}\times S,\gamma^*V)$ and 
\beq
i_0\colon S\times \R\hookrightarrow S\times \R^{1|1}\label{eq:i0}
\eeq
 the canonical inclusion of the fiberwise reduced manifold. These are a vector bundle version of the standard Taylor expansion of a function on $\R^{1|1}$ in the odd variable~$\theta$, and it's easy to check (e.g., in a local trivialization) that a section is determined by its components. Note that $s_0$ and $s_1$ are sections of $V$ and $\Pi V$ pulled back along a \emph{super} family of \emph{ordinary} paths, $i_0\circ \gamma\colon S\times \R\to M$. 

\begin{prop} \label{prop:parallelcomp}
A section is super parallel if and only if its components satisfy
$$
s_1=0,\quad (i_0^*\gamma^*\nabla)_{\partial t}(s_0)=-\frac{1}{2} i_0^*(\iota_D\iota_D(\gamma^*F))s_0
$$
where $F\in \Omega^2(M;\End(V))$ is the curvature 2-form of $\nabla$ and $\iota_D$ is contraction with~$D$. 
\end{prop}

\bp
To simplify notation below, set $\nabla=\gamma^*\nabla$. Using the definitions~\eqref{eq:components}, we compute
$$
(\nabla_Ds)_0=s_1,\quad (\nabla_Ds)_1=i_0^*(\nabla_D(\nabla_Ds))=i_0^*((\nabla_t+\frac{1}{2}(\gamma^*F)(D,D))s)=(\nabla_t+\frac{1}{2}\gamma^*F(D,D))i_0^*s
$$
where we used
$$
\nabla_D\nabla_D+\nabla_D\nabla_D-\nabla_{[D,D]}=(\gamma^*F)(D,D), \quad [D,D]=2\partial_t.
$$
So the equation $\nabla_Ds=0$ is equivalent to the claimed equations in the proposition. 
\ep

Since $s_1\equiv 0$ for all $t$, the above lemma allows one to view super parallel transport as a single differential equation involving $s_0$ over $\R\times S$. In local coordinates, one can check that this reduces to an ordinary first order differential equation so the usual existence and uniqueness of solutions apply for an initial condition $s_0(t_0)$. 

\subsection{Super holonomy on constant super loops and the Chern character}

Now we study the super parallel transport formula on the constant super loops $\SM(\R^{0|1},M)\subset \sLoop(M)\subset\sPath(M)$. Let $p\colon \Pi TM\to M$ be the projection, and note that after identifying $\Pi TM\cong \SM(\R^{0|1},M)$, we have $i_0\circ \ev=p$ for $\ev$ and $i_0$ as in \eqref{eq:ev} and \eqref{eq:i0}, respectively.

\begin{lem} Let $(V,\nabla)$ be a vector bundle with connection on~$M$. Under the identification,
$$
\Gamma(\SM(\R^{0|1},M),p^*\End(V))\cong \Gamma(\Pi TM,p^*\End(V)\cong \Omega^\bullet(M,\End(V)),
$$
the section $i_0^*(\ev^*F(D,D))\in \Gamma(\SM(\R^{0|1},M),p^*\End(V))$ is $-2F\in \Omega^\bullet(M,\End(V))$. 
\end{lem}

\bp It suffices to check the statement locally. Below, for $f\in C^\infty(M)$ let $\delta f\in C^\infty(\Pi TM)$ denote the function associated to the 1-form $df\in \Omega^1(M)$. Hence, $\ev^*f=f+\theta \delta f$ in this notation. Pulling back the 2-form $F=dfdg\in \Omega^2(M)$ along~$\ev$, we compute $i_0^*(\iota_D\iota_D\ev^*F)$
\beq
i_0^*(\iota_D\iota_Dd(f+\theta \delta f)d(g+\theta \delta g))&=&i_0^*(\iota_D\iota_D (df+d\theta \delta f+\theta d\delta f)(dg+d\theta \delta g+\theta d\delta g))\nonumber \\
&=&-2\delta f\delta g\in C^\infty(\Pi TM)\cong C^\infty(\SM(\R^{0|1},M))\nonumber
\eeq
Working in a local trivialization of~$V$, this shows that~$i_0^*(\ev^*F(D,D))$ is identified with the claimed endomorphism-valued form. 
\ep

From the lemma, the super parallel transport equation on constant super paths reads 
$$
s_1=0,\quad \frac{ds_0}{dt}=F s_0,\quad s_1\in \Omega^\bullet(M,p^*\Pi V), \ s_0\in \Omega^\bullet(M,p^*V)
$$
where now~$F$ does not depend on $t$. We have the solution,
$$
s_1(t)=0, \quad s_0(t)=\exp(tF)s_0(0),
$$
for an initial value $s_0(0)$. Hence in this case super parallel transport is identified with
\beq
{\rm sPar}(V,\nabla,t)=\exp(tF)\in \Omega^\bullet(M,\End(V))\cong \Gamma(\SM(\R^{0|1},M),\End(p^*V))\label{eq:susyhol}
\eeq
so that ${\rm sPar}(V,\nabla,t)(s_0(0))=s_0(t)$ solves the initial value problem. 

\begin{defn} The \emph{super holonomy} on the constant super loops is the endomorphism \eqref{eq:susyhol} at $t=1$
$$
\shol(V,\nabla):={\rm sPar}(V,\nabla)(1)\in C^\infty(\SM(\R^{0|1},M))\cong \Omega^\bullet(M).
$$
 \end{defn}

From~\eqref{eq:susyhol}, the following is clear. 

\begin{thm} 
The trace of the super holonomy on constant super loops coincides with the Chern character, 
$${\rm Tr}(\shol(V,\nabla))={\rm Tr}(e^{F})\in\Omega^{\ev}_{\rm cl}(M)\subset C^\infty(\SM(\R^{0|1},M)),$$ 
as a closed even form on~$M$. \end{thm}

\section{Equivariant super paths and equivariant super parallel transport}

For a $G$-manifold $M$, the natural generalization of a super path in $M$ is maps into the quotient stack, 
\beq
(P,\gamma_G)\colon \R^{1|1}\times S\to M\sq G\quad \iff\quad \R^{1|1}\times S\leftarrow P\stackrel{\gamma_G}{\to} M\label{eq:maptostack}
\eeq
where $P\to \R^{1|1}\times S$ is a $G$-bundle and $P\to M$ is a $G$-equivariant map. We will make sense out of super parallel transport for sections $s\in \Gamma(P,\gamma_G^*V)$ for $(V,\nabla)$ an equivariant vector bundle with invariant connection. For this to work, we require a lift of the vector field~$D$ to~$P$. A connection~$A$ on~$P$ gives a unique such lift of~$D$ which we denote by $D_A$. Since the connection is crucial part of the story, we define an $S$-family of \emph{equivariant super paths} to be a map as in~\eqref{eq:maptostack} together with a connection $A$ on the $G$-principal bundle~$P$. 

\begin{defn} \label{defn:equivparallel}
For an equivariant super path, let $D_A$ denote the horizontal lift of $D$ using the connection~$A$. A section is \emph{equivariantly super parallel} if
\beq
\nabla_{D_A} s=0, \quad s\in \Gamma(P,\gamma^*V),\label{eq:anequivparal}
\eeq
for $\gamma\colon P\to M$ a smooth map. We call~\eqref{eq:anequivparal} \emph{the equivariant super parallel transport equation}. 
\end{defn}

As per the ordinary Chern character, we will focus attention on \emph{constant} equivariant paths with \emph{constant} connections. Constancy will mean that $\R^{1|1}\times S\to M\sq G$ factors through the quotient $(S\times \R^{1|1}/\E)$ (so $P$ pulls back from a $G$-bundle on $\R^{0|1}\times S$) and the connection $A$ is gauge-equivalent to an $S$-family of connections on $G$-bundles on~$\R^{0|1}$. 

\subsection{The equivariant super path space}\label{sec:eqpath}

For a trivial $G$-bundle $G\times \Sigma\times S\to \Sigma\times S$, an \emph{$S$-family of connections} can be identified with a 1-form,
$$
A\in \Omega^1_S(\Sigma;\mathfrak{g}):= \Omega^1(\Sigma;\mathfrak{g})\otimes C^\infty(S).
$$
In the case at hand with $\Sigma=\R^{1|1}$, all $G$-bundles are trivializable (locally in $S$) so a map
$$
\R^{1|1}\times S\to M\sq G
$$
into the quotient stack is isomorphic (locally in $S$) to an equivariant map
$$
\gamma_G\colon G\times \R^{1|1}\times S\to M.
$$
With this simplification in play, we make the following definition. 

\begin{defn} \label{defn:equivpath} For $M$ a manifold with $G$-action, the presheaf of \emph{equivariant super paths}, denoted $\sPath(M\sq G)$, has as $S$-points pairs $(\gamma_G,A)$ where $\gamma_G\colon G\times \R^{1|1}\times S\to M$ is an equivariant map and $A$ is an $S$-family of connections~$A\in \Omega^1_S(\R^{1|1};\mathfrak{g})$. 
\end{defn}

\begin{rmk}\label{rmk:gauge}
In truth, equivariant super paths form a sheaf of groupoids (i.e., a stack) where isomorphisms come from \emph{gauge transformations} between equivariant super paths. These are maps $g\colon \R^{1|1}\times S \to G$ that acts on the map $\gamma_G$ by precomposing with the $G$-action on the fibers of the $G$-bundle, and changes $A$ by the gauge transformation formula,
\beq
A\mapsto \Ad_gA+g^*\Theta\label{eq:gaugetransf}
\eeq
where $\Theta$ is the Maurer--Cartan form; $g^*\Theta=d\log(g)$ for matrix groups. Gauge transformations fit together with the translational action $\E^{1|1}$ on $\R^{1|1}\times S$ to give an action of $\E^{1|1}\ltimes \SM(\R^{1|1},G),$ on $\sPath(M\sq G)$. We'll make limited use of this additional structure, so for simplicity work with the presheaf (of sets) $\sPath(M\sq G)$. However, working with the stack has its benefits. For example, the groupoid of $G$-bundles with connection over~$\R^{0|1}$ provides a model for equivariant de~Rham cohomology; see~\cite{HSST} for details. 
\end{rmk}

Definition~\ref{defn:equivpath} gives an isomorphism of presheaves, $\sPath(M\sq G)\cong \sPath(M)\times \Omega^1(\R^{1|1};\mathfrak{g})$. Indeed, any map $\gamma_G$ necessarily factors as 
\beq
G\times \R^{1|1}\times S\stackrel{\id_G\times \gamma}{\longrightarrow} G\times M\stackrel{\mu}{\longrightarrow} M\label{eq:stack}
\eeq
so is determined by a map $\gamma\colon \R^{1|1}\times S\to M.$

\subsection{Equivariant super parallel transport for constant $G$-connections}

An $\E$-invariant connection pulls back from $\R^{0|1}$, meaning it is in the image, $A\in \Omega^1_S(\R^{0|1};\mathfrak{g})\to \Omega^1_S(\R^{1|1};\mathfrak{g})$. For a choice of coordinate $\theta$ on $\R^{0|1}$ such connections take the form,
$$
A=d\theta\otimes \alpha +\theta d\theta\otimes a,\quad \alpha\in \Pi\mathfrak{g}(S),\ a\in \mathfrak{g}(S). 
$$
As per Remark~\ref{rmk:gauge}, we can see that the gauge transformation $e^{-\theta\alpha}\colon S\times \R^{1|1}\to G$ has the effect
\beq
A\mapsto A'=d\theta\otimes \alpha+\theta d\theta\otimes a-d\theta\otimes \alpha=\theta d\theta\otimes a.\label{eq:gaugetransf}
\eeq
Hence, any $G$-connection on $\R^{0|1}$ is gauge equivalent to one of the form $A=\theta d\theta\otimes a$. This justifies the following definition. 

\begin{defn}
A \emph{constant $G$-connection} for an equivariant super path is one of the form $A=\theta d\theta\otimes a$ for $a\in \mathfrak{g}(S)$. 
\end{defn}

The horizontal lift $D_A$ of $D$ with respect to the connection $A=\theta d\theta\otimes a$ is
$$
D_A:=D+\iota_DA=D-\theta a
$$
where we regard $a\in\mathfrak{g}(S)$ as determining a vector field on $G\times S$, and hence $D_A$ determines a vector field on $G\times \R^{1|1}\times S$. 
\begin{rmk} Another interpretation of~$D_A$ is that it is the essentially unique lift of~$D$ preserving a constant $G$-connection~$A$. Indeed, the action of $(s,\eta)\in \E^{1|1}(S)$ on $A$ generated by $D$ is
$$A\mapsto (\theta+\eta)d\theta \otimes a=A+\eta d\theta \otimes a,$$
while the action of gauge transformations by $\exp(-\eta \theta a) \in \Pi \mathfrak{g}(\E^{1|1}\times \R^{1|1}\times S)$ is 
$$ A\mapsto A+d\log(\exp(-\eta \theta a))=A-\eta d\theta\otimes a.$$
Hence, the sums of the derivatives of these actions at $(s,\eta)=(0,0)$ preserves~$A$. The derivative of the first action at zero is the vector field $D$, and the derivative of the second action at zero is $-\theta a$, showing that $A$ is invariant under $D_A=D-\theta a$. This vector field is closely related to the Cartan differential in a gauged field theory model for equivariant de~Rham cohomology; see~\cite{HST}. 
\end{rmk} 

We now analyze the equivariant super parallel transport for constant $G$-connections in component fields. For $i_0\colon S\times \R\times G\hookrightarrow \R^{1|1}\times S\times G$, we define
$$
s_0:=i_0^*s,\quad s_1:=i_0^*((\gamma^*\nabla)_{D_A} s).
$$

\begin{prop} \label{prop:equivparallelcomp}
For a constant $G$-connection~$A$, a section is equivariantly super parallel if and only if its components satisfy
$$
s_1=0,\quad (\gamma^*\nabla)_{\partial t}(s_0)=-\left(\frac{1}{2} i_0^*((\gamma^*F)(D,D))-\iota_a (i_0^*\gamma^*\nabla)\right) s_0
$$
where $F$ is the curvature 2-form of $\nabla$, $a$ is the (even) vector field on $G\times \R\times S$ associated with the connection~$A$, and $\partial_t$ is the vector field on $G\times \R\times S$ associated with translations in $\R$. Hence, this coincides with the non-equivariant super parallel transport equation with an extra term involving~$a$. 
\end{prop}

\bp
Applying the same argument as in the non-equivariant setting, we find that the super parallel transport equation is equivalent to 
$$
s_1=0,\quad i_0^*((\gamma^*\nabla_{\partial_t-a}+\frac{1}{2}(\gamma^*F)(D_A,D_A))s_0)=0
$$
where we use that $[D_A,D_A]=2(\partial_t-a)$. Since $\partial_t-a$ pushes forward from $S\times \R\times G$, we have
$$
i_0^*(\gamma^*\nabla_{\partial_t-a}s_0)=(i_0^*\gamma^*\nabla)_{\partial_t-a}(s_0). 
$$
With $(i_0^*\gamma^*\nabla)=:\nabla^0$, we can re-express this as
$$
\nabla_{\partial_t}^0 (s_0)=-\frac{1}{2}i_0^*(\iota_{D_A}\iota_{D_A} \gamma^*F)s_0+\nabla^0_a s_0=-\frac{1}{2}i_0^*(\iota_{D_A}\iota_{D_A} \gamma^*F)s_0+\iota_a \nabla^0 s_0.
$$
Now, 
$$
i_0^*(\iota_{D_A}\iota_{D_A} \gamma^*F)=i_0^*((\iota_D+\theta\iota_a)(\iota_D+\theta\iota_a)\gamma^*F)=i_0^*(\iota_D\iota_D\gamma^*F)
$$
where we use that $i_0^*(\theta)=0$, and so the terms involving $\theta\iota_a$ can be discarded. This gives the claimed equation. 
\ep

\subsection{Equivariant super parallel transport on constant equivariant super paths}

\begin{defn} The presheaf of \emph{constant equivariant super paths}, denoted $\sPath_0(M\sq G)$, has as $S$-points constant $G$-connections $A=\theta d\theta\otimes a\in \in \Omega^1_S(\R^{1|1}/\E;\mathfrak{g})$ and equivariant maps 
$$
\gamma_G\colon G\times \R^{1|1}\times S\twoheadrightarrow (G\times \R^{1|1}\times S)/\E\cong G\times \R^{0|1}\times S \to M.
$$
\end{defn}

Above, the map $\gamma_G$ is determined by $\gamma\colon S\times \R^{0|1}\to M$
$$
\gamma_G\colon G\times \R^{1|1}\times S\twoheadrightarrow G\times \R^{0|1}\times S\stackrel{\id_G\times \gamma}{\longrightarrow} G\times M\stackrel{\mu}{\to }M.
$$
This gives a description of constant equivariant super paths as the presheaf
$$
\mathfrak{g}\times \SM(\R^{0|1},M)\cong \sPath_0(M\sq G)\hookrightarrow \sPath(M\sq G).
$$
The universal map $\gamma_G$ to~$M$ is determined by the evaluation $\gamma=\ev$ 
$$
\gamma_G\colon G\times \R^{1|1}\times (\mathfrak{g}\times \SM(\R^{0|1},M))\stackrel{{\rm pr}}{\to} G\times \R^{0|1}\times \SM(\R^{0|1},M)\stackrel{\ev}{\to} G\times M\stackrel{\mu}{\to} M,
$$
where ${\rm pr}$ projects out $\R^1$ and $\mathfrak{g}$. 

\begin{prop}\label{prop:ODE} 
For $(V,\nabla)$ an equivariant vector bundle with invariant connection on~$M$, a section over the universal family of constant equivariant super paths is equivariantly super parallel if and only if its first component vanishes, $s_1\equiv 0$, and its zero component
satisfies 
$$
(\nabla_{\partial_t}s_0)(t)=(F+[\iota_a,\nabla])s_0(t)
$$
where for each $a\in \mathfrak{g}$, $F+[\iota_a,\nabla]$ is a differential forms on $G\times M$ valued in endomorphisms of $\mu^*V$ and we have identified 
\beq
&&s_0(t)\in C^\infty(\mathfrak{g})\otimes \Omega^\bullet(M)\otimes_{C^\infty(M)} \Gamma(G\times M;\mu^*V)\cong \Gamma(G\times \mathfrak{g}\times \SM(\R^{0|1},M),\gamma_G^*V)\label{eq:translate}
\eeq
using the isomorphism functions on $\SM(\R^{0|1},M)$ and differential forms on $M$. 
\end{prop}

\bp From Proposition~\ref{prop:equivparallelcomp}, we need only translate the equation involving~$s_0$ into differential forms via~\eqref{eq:translate}. Since the connection (and therefore its curvature) are $G$-invariant, the pullback of $F$ along the action map $\mu\colon G\times M\to M$ coincides with the pullback along the projection to~$M$. This allows us to identify $F$ with the pullback of the curvature term from the nonequivariant case along the projection $G\times \SM(\R^{0|1},M)\to \SM(\R^{0|1},M).$

For the remaining term, we write a section of the pullback of~$V$ along $\gamma_G$ as $f\otimes s$ for $s$ a section of $\mu^*V$ over $G\times M$ and $f\in C^\infty(\Pi TM)$. Then we have
\beq
\iota_a \nabla (f\otimes s)=\iota_a(df\otimes s+f\otimes \nabla s)=\mathcal{L}_a f\otimes s+f\otimes\nabla_a s.\label{eq:Lie}
\eeq
The Lie derivative term is the Cartan formula together with the fact that $f$ has cohomological degree zero. Note there are no signs because $f$ has cohomological degree zero and $\nabla$ has super degree zero. 

Viewing $C^\infty(\Pi TM)$ as differential forms on $M$ and a section of the pullback of $V$ as a differential form $\omega$ valued in a section $s$ on $\R\times M$, we have
$$
\iota_a\nabla (\omega\otimes s)=\iota_a(d\omega\otimes s\pm \omega\otimes \nabla s)=\iota_a d\omega\otimes s + \omega\otimes \nabla_a s\pm \iota_a\omega \otimes \nabla s.
$$
But then $\iota_a d=\mathcal{L}_a-d\iota_a$, and since the cohomological degree of $\omega$ need not be zero, $d\iota_a$ and $\nabla \iota_a$ need not act by zero:
$$
(\nabla \iota_a)(\omega\otimes s)=\nabla(\iota_a \omega \otimes s)=d\iota_a \omega \otimes s\mp \iota_a \omega \otimes \nabla s. 
$$
This gives us 
$$
[\iota_a,\nabla](\omega\otimes s)=\mathcal{L}_a\omega\otimes s+\omega\otimes \nabla_a s.
$$
Comparing with~\eqref{eq:Lie}, this completes the proof. 
\ep

\section{Constant equivariant super loops and the equivariant Chern character}

An \emph{equivariant loop} is a path $\gamma\colon \R\to M$ and a group element $g\in G$ such that 
\beq
\gamma(t)=g\cdot \gamma(t+1). \label{eq:eqloop}
\eeq
Given an equivariant vector bundle $V\to M$ and a linear map $V_{\gamma(0)}\to V_{\gamma(1)}$ between fibers from parallel transport, we get an endomorphism of $V_{\gamma(0)}$ by composing with the map~$g\colon V_{\gamma(1)}\to V_{\gamma(0)}$ gotten from the equivariant structure. This gives a version of \emph{equivariant holonomy}. This seems to be a well-known object, though we have been unable to pinpoint its origin.

We can repackage an equivariant loop $\gamma$ as above in terms of maps $S^1\to M\sq G$: this has as data a principal $G$-bundle $P\to S^1$ with a $G$-equivariant map~$\phi\colon P\to M$. If we equip $P$ with a connection and choice of basepoint $p\in P$, there is a unique lift of the standard covering map $\R\to \R/\Z\cong S^1$ to a path $\R\to P$. Denote the composition $\R\to P\to M$ by $\gamma$. Since the lifts of $1\in \R$ and $0\in \R$ are in the same fiber of $P$ and $\phi\colon P\to M$ is $G$-equivariant, there is a unique element $g\in G$ such that $g\cdot \gamma(1)=\gamma(0)$. This gives the equivariant path of the previous paragraph one which the equivariant holonomy is defined. 

Below we construct a super equivariant analog of the construction above. The output is a function on constant equivariant super loops, which we identify with a subsheaf
$$
\sLoop_0(M\sq G)\subset G\times \mathfrak{g}\times \SM(\R^{0|1},M).
$$
We identify this function on $\sLoop_0(M\sq G)$ with the bouquet of Chern characters.

\subsection{Constant equivariant super loops}

\begin{defn} An $S$-family of \emph{equivariant super loops} in $M$ is $(P,A,\phi)$ where $P\to S\times \R^{1|1}/\Z$ is a $G$-principal bundle, $A$ is an $S$-family of connections on~$P$, and $\phi\colon P\to M$ is an equivariant map. \end{defn}

\begin{lem} \label{lem:equivloops}
An $S$-family of equivariant super paths $(\gamma_G,A)$ and an $S$-point $h\in G(S)$ determine an $S$-family of equivariant super loops if $A$ is invariant for the action~\eqref{eq:gaugetransf} of gauge transformations by~$h$ and $\gamma_G$ is invariant under the $\Z$-action determined by the standard $\Z<\E<\E^{1|1}$-action on $\R^{1|1}$ and the $\Z$-action on $G\times S$ generated by $h$
\beq
\begin{tikzpicture}[baseline=(basepoint)];
\node (B) at (5,0) {$\Z\times G\times \R^{1|1}\times S$};
\node (D) at (5,-1.5) {$G\times \R^{1|1}\times S$};
\node (E) at (9,-.75) {$M.$};
\draw[->] (B) to node [left=2pt] {${\rm act}$} (D);
\draw[->] (B) to node [above=2pt] {$\gamma_G\circ{\rm pr}$} (E);
\draw[->] (D) to node [below=2pt] {$\gamma_G$} (E);
\end{tikzpicture}
\nonumber
\eeq
\end{lem}
\bp In this case, $h$ determines an $S$-family of $G$-bundles $P=(G\times \R^{1|1}\times S)/\Z\to \R^{1|1}/\Z\times S$ with the property that $A$ descends to a connection on~$P$ and $\gamma_G$ descends to an equivariant map $P\to M$ on the quotient. \ep

Unlike the nonequivariant case, constant equivariant super paths are not automatically equivariant super loops. There are two reasons for this: (1) an equivariant loop is \emph{more data} than an equivariant path (we need to specify $h\in G$), and (2) there is a nontrivial condition to promote a constant equivariant super path to a constant equivariant super loop. 

\begin{defn} The presheaf of \emph{constant equivariant super loops} is the sub presheaf 
$$
\sLoop_0(M\sq G)\subset G\times \mathfrak{g} \times \SM(\R^{0|1},M)\
$$
for which the $S$-point of $\mathfrak{g}\times \SM(\R^{0|1},M)\cong \sPath_0(M\sq G)$ together with the $S$-point of $G$ define an $S$-family of equivariant super loops in $M$ as in Lemma~\ref{lem:equivloops}. 
\end{defn}

Explicitly, this subsheaf is characterized as having $S$-points $h\in G(S)$, $a\in \mathfrak{g}(S)$ and $\R^{0|1}\times S\to M$ satisfing
\begin{enumerate}
\item $\Ad(h,a)=(h,a)$ where $\Ad\colon G\times \mathfrak{g}\to G\times \mathfrak{g}$ is the Cartesian product of the projection $G\times \mathfrak{g}\to G$ and the adjoint action $G\times \mathfrak{g}\to \mathfrak{g}$;
\item there is a factorization $\R^{0|1}\times S\to M^h\subset M$ through the fixed point set of~$h$. 
\end{enumerate}

There is a projection map $\sLoop_0(M\sq G)\to G$, and we think of $\sLoop_0(M\sq G)$ as a sort of bundle over~$G$. The total space is typically not representable, as the fixed point sets $M^h$ need not vary continuously with $h\in G$. However, for each $h\in G$ we get a representable fiber
$$
\sLoop_0(M\sq G)_h:=\mathcal{Z}(h)\times \SM(\R^{0|1},M^h)\hookrightarrow \sLoop_0(M\sq G)
$$
where $\mathcal{Z}(h)\subset \mathfrak{g}$ consists of the $\Ad_h$-invariant subspace of $\mathfrak{g}$. 

\subsection{Super holonomy on constant equivariant super loops}

Let $(P,A,\phi)$ be an equivariant super loop in~$M$. A basepoint $p\in P$ determines a lifting $\exp(D_A)$ of the covering map $S\times \R^{1|1} \to S\times \R^{1|1}/\Z$ to~$P$ using the connection,\footnote{Explicitly, $\exp(D_A)$ sends $\R^{1|1}\times S$ to the $\E^{1|1}$-orbit of the basepoint for the action generated by $D_A$. } and we have the composition
\beq
\gamma_0\colon \R\times S\stackrel{i_0}{\hookrightarrow} \R^{1|1}\times S \stackrel{\exp(D_A)}{\to} P\to M.\label{eq:gamma0}
\eeq
Let $V_0$ denote the fiber of $\gamma_0^*V$ at $\{0\}\times S$ and $V_1$ denote the fiber at $\{1\}\times S$. By Proposition~\ref{prop:equivparallelcomp}, the restriction of equivariantly super parallel sections to the source $\R\times S$ of~\eqref{eq:gamma0} satisfy an initial value problem. Hence they determine a linear map~$V_0\to V_1$ over~$S$ sending initial conditions at $t=0$ to their solutions at~$t=1$. Moreover, the image of $\{0\}\times S$ and $\{1\}\times S$ in~$P$ are in the same fiber, so they differ by a unique element~$g\in G(S)$. Using the equivariant structure on~$V$ and that $\phi$ is a $G$-equivariant map, $g\in G(S)$ defines a linear map $V_1\to V_0$ over $S$. 

\begin{defn} 
The \emph{equivariant super holonomy} is the endomorphism of~$V_0$ gotten by the composition of the map $V_0\to V_1$ determined by equivariantly super parallel sections, and the map $V_1\to V_0$ determined by the unique~$g\in G(S)$ in the notation above. 
\end{defn}

\begin{lem} \label{lem:superholg} For an $S$-family of constant equivariant super loops determined by $h\in G(S)$ with connection associated with $a\in \mathfrak{g}(S)$, the group element $g\in G(S)$ in the definition of super holonomy is $h\cdot e^{a}\in G(S)$. 
\end{lem}
\bp 
We start by considering an $S$-family of constant equivariant super paths, $S\to \sPath_0(M\sq G)\cong \mathfrak{g}\times \SM(\R^{0|1},M)$. The vector field~$D_A$ and the basepoint of $G\times \R^{1|1}\times \sPath_0(M\sq G)$ associated to $e\in G$ and $(0,0)\in \R^{1|1}$ determines 
\beq
\gamma_0\colon \R\times S\stackrel{i_0}{\hookrightarrow} \R^{1|1}\times S\stackrel{D_A}{\to} G\times \R^{1|1}\times S\stackrel{\gamma_G}{\to} M\label{eq:superlift}
\eeq
where $\exp(D_A)$ is a super path lifting.  Post-composing $\gamma_0$ with the projection to~$\mathfrak{g}\times G$, the image of $\{1\}\times S$ differs from the image of $\{0\}\times S$ by $e^{a}\in G$, for~$a\in\mathfrak{g}$. In particular, the action of $e^a$ on the image of $\{1\}\times S$ in $M$ is the image of $\{0\}\times S$ under $\gamma_0$, where we have used the trivialization of the $G$-bundle to identify fibers. 

Given an equivariant super loop, we have the additional data~$h\in G(S)$ that generates a $\Z$-action with quotient $P=(G\times \R^{1|1}\times S)/\Z$. This $h\in G(S)$ is precisely an identification between the fibers over $\{0\}\times S\subset \R^{1|1}\times S$ and $\{1\}\times S\subset \R^{1|1}\times S$. Taken together with case where we use the trivialization to identify fibers (i.e., when $h=e$), the unique group element identifying the images of $\{0\}\times S$ and $\{1\}\times S$ is the product, $h\cdot e^{a}$, proving the lemma. 
\ep

\begin{proof}[Proof of Theorem~\ref{thm1}]
Before taking a trace, we construct a linear map related to the super holonomy on constant equivariant paths equipped with a choice of~$g\in G$, i.e., on
$$
G\times \mathfrak{g}\times \SM(\R^{0|1},M)\cong G\times \sPath_0(M\sq G) \supset \sLoop_0(M\sq G).
$$ 
This yields a map between different bundles that when restricted to the constant equivariant super loops gives a section of an endomorphism bundle. 

In the notation of \eqref{eq:superlift}, we consider $\gamma_G\circ \gamma_0$ for $S=\sPath(M\sq G)$,
$$
\R\times \sPath_0(M\sq G) \hookrightarrow \R^{1|1}\times \sPath_0(M\sq G)\stackrel{\exp(D_A)}{\longrightarrow} G\times \R^{1|1}\times \sPath_0(M\sq G)\stackrel{\gamma_G}{\to} M.
$$
We can identify this composition with
$$
\R\times \sPath_0(M\sq G)\cong (\R\times\mathfrak{g})\times \Pi TM\stackrel{\exp\times {\rm pr}}{\longrightarrow} G\times M\stackrel{\mu}{\to} M
$$
where $\exp\colon \R\times \mathfrak{g}\to G$ is the exponential map $(t,a)\mapsto e^{-ta}$, and ${\rm pr}\colon \Pi TM\to M$ is the projection. From Proposition~\ref{prop:ODE}, the initial value problem for equivariant super parallel transport has solutions 
$$
s_0(t)=e^{t(F+[\iota_a,\nabla])}s_0(0),\quad s_0\in \Gamma((\gamma_G\circ \gamma_0)^*V)
$$
for the initial condition~$s_0(0)$. As $t$ varies, this gives maps between the pullback of $V$ to $\{0\}\times S$ and the pullback of this map post-composed with the action of $e^{-ta}$ on $M$. When $t=1$, we get the map 
$$
e^{F+[\iota_a,\nabla]}\colon V_0\to V_1.
$$
From the definition of equivariant super holonomy together with Lemma~\ref{lem:superholg}, we postcompose the above with the action~$\mu^V_{h\cdot e^{ta}}\colon V_1\to V_0$. Letting~$h$ vary, this gives
$$
\mu^V e^{(F+[\iota_a,\nabla]+\mathcal{L}^V_a)}\in \Gamma(G\times \mathfrak{g}\times \Pi TM,\Hom(p^*V,\mu^*V))\cong C^\infty(G\times \mathfrak{g}\times M;\mu^*V)\otimes_{C^\infty(M)}\Omega^\bullet(M),
$$
where $\mathcal{L}_a^V$ generates the $a$-action on~$V$. This gives the section in the statement of Theorem~\ref{thm1}. 

Now we restrict the above to equivariant super loops for a fixed $h\in G$. This requires~$\Ad_h(a)=a$ and $\SM(\R^{0|1},M^h)\subset \SM(\R^{0|1},M)$. From the definition of equivariant super holonomy, we obtain
$$
\mu^V(h)\circ e^{(F+[\iota_a,\nabla]+\mathcal{L}^V_a)} \in \Gamma(\Pi TM^h,\End(p^*V))\cong \Omega^\bullet(M^h;\End(V)). 
$$
Here it is important that the action $\mu^V(h)$ is through endomorphisms of the restriction of~$V$ to~$M^h$. The traces of the endomorphisms above are precisely pieces of the bouquet of differential forms comprising the bouquet of Chern characters of $(V,\mu^V,\nabla)$; see~\ref{eg:bouquet}.
\ep

\begin{rmk} If we consider constant equivariant super loops of circumference $\epsilon$ rather than circumference~1 and we set $h=e$, then the super holonomy is 
$$
e^{\epsilon (F+[\iota_a,\nabla]+\mathcal{L}^V_a)}=\id+\epsilon (F+[\iota_a,\nabla]+\mathcal{L}^V_a)+O(\epsilon^2).
$$
The coefficient of $\epsilon$ above is precisely the equivariant curvature. In this sense, equivariant super parallel transport around infinitesimal loops is the equivariant curvature. 
\end{rmk}

\appendix

\section{The bouquet of Chern characters}\label{appen}

Following~\cite{BGV}, for $(V,\nabla)$ an equivariant vector bundle with invariant super connection, the \emph{equivariant curvature} is a differential form-valued function on $\mathfrak{g}$ whose value at $X\in\mathfrak{g}$ is
$$
F(X):=(\nabla-\iota_{X})^2+\mathcal{L}^V(X)
$$
where $\mathcal{L}^V(X)$ denotes the infinitesimal $G$-action on sections of~$V$ by~$X$. The delocalized Chern character of $(V,\nabla)$, as studied by Block and Getzler~\cite{BlockGetzler} and Vergne~\cite{Vergne}, takes values in a somewhat sophisticated version of equivariant de~Rham cohomology that we review now.

\begin{defn} A \emph{bouquet of equivariant differential forms} on $M$ assigns to each $g\in G$ a closed equivariant differential form $\alpha_g\in \Omega_{\rm cl}^\bullet (M^g,\mathcal{Z}(g))$ satisfying
\begin{enumerate} 
\item $h\cdot \alpha_g=\alpha_{hgh^{-1}}$ for all $h\in G$ and
\item for sufficiently small $\epsilon\in \R$, for all $X,Y\in \mathcal{Z}(g)$ we have an equality
$$
\alpha_{ge^{\epsilon X}}(Y)=\alpha(\epsilon X+Y)|_{M^{ge^{\epsilon X}}}
$$
 of differential forms  on $M^{ge^{\epsilon X}}$. 
\end{enumerate}
\end{defn}

\begin{ex} For a $G$-representation $\rho\colon G\to \End(V)$, 
$$
\alpha_g(X):={\rm Tr}(\rho(g e^{X}))
$$
defines a bouquet of equivariant differential forms on $\pt$ \cite{BlockGetzler,Vergne}. Note that on restriction to $X=0\in \mathcal{Z}(g)$, this is the character of the representation. 
\end{ex}

\begin{ex} \label{eg:bouquet}
For $(V,\nabla)$ an equivariant vector bundle with invariant super connection, 
\beq
{\rm ch}(V,\nabla)_g(X)={\rm Tr}(\mu^V(g)\circ e^{F(X)|_{M^g}})\label{eq:globalCC}
\eeq
is a bouquet of equivariant differential forms on~$M$ \cite{BlockGetzler,Vergne}.
\end{ex}

\bibliographystyle{amsalpha}

\begin{thebibliography}{EKMM}

\bibitem[AG83]{Alvarez} Alvarez-Gaum\'e, {\em Supersymmetry and the Atiyah-Singer index theorem,} Communications in Mathematical Physics 90 (1983), 161-173.

\bibitem[AS69]{AS} M. Atiyah and G. Segal, {\em Equivariant K-theory and completion.} Journal of Differential Geometry (1969). 

\bibitem[B85]{B85} J-M.~Bismut,
{\em Index theorem and equivariant cohomology on the loop space}, Comm. Math. Phys. {\bf 98} (1985), no. 2, 213-237.


\bibitem[BG94]{BlockGetzler} J. Block and E. Getzler, {\em Equivariant cyclic homology and equivariant differential forms}, Ann. Sci. Ecole Normale Sup. 4 (1994).

\bibitem[BGV92]{BGV} N. Berline, E. Getzler, and M. Vergne, {\em Heat kernels and Dirac operators}, Springer, 1992.\

\bibitem[DM99]{DM} P. Deligne and J. Morgan, {\em Notes on supermanifolds}, Quantum Fields and Strings: A Course for Mathematicians, Volume 1 (P. Deligne, P. Etingof, D. Freed, L. Jeffrey, D. Kazhdan, J. Morgan, D. Morrison, and E. Witten, eds.), American Mathematical Society, 1999.

\bibitem[DV93]{DV} M. Duflo and M. Vergne, {\em Cohomologie \'{e}quivariante et descente}, Ast\'{e}risque, 215, 5-108 (1993).

\bibitem[D06]{Florin} F. Dumitrescu, {\em Superconnections and parallel transport}, Ph.D. Thesis, 2006. 


\bibitem[FHT07]{FHT} D. Freed, M. Hopkins and C. Teleman, {\em Twisted equivariant K-theory with complex coefficients}, Topology, 1 (2007). 

\bibitem[Han08]{Han} F. Han, {\em Supersymmetric QFTs, SuperLoop Spaces and Bismut--Chern Character.} PhD thesis 2008.

\bibitem[HSST]{HSST}  F. Han, C. Schommer-Pries, S. Stolz and P. Teichner, {\em Equivariant cohomology from
gauged field theories}, to appear. 

\bibitem[HST]{HST}  F. Han, S. Stolz and P. Teichner, {\em The Bismut--Chern character form as dimensional reduction}, to appear. 

\bibitem[Ro03]{Rosu} I. Rosu, {\em Equivariant K-theory and equivariant cohomology.} With an appendix by Knutson
and Rosu. Math. Z. 243 (2003), 423-448.

\bibitem[ST11]{ST11} S. Stolz and P. Teichner, {\em Supersymmetric field theories and generalized cohomology}, Mathematical
Foundations of Quantum Field and Perturbative String Theory (B. Jur\v co, H. Sati, U. Schreiber,
ed.), Proceedings of Symposia in Pure Mathematics, 2011.

\bibitem[Ver94]{Vergne} M. Vergne, Geometric quantization and equivariant cohomology, ECM: Proceedings of the First
European Congress of Mathematics (A. Joseph, F. Mignot, F. Murat, B. Prum, and R. Rentschler, eds.), Progress in Mathematics, vol. 119, Birkh\"auser, 1994.


\end{thebibliography}


\end{document}